\documentclass[12pt]{amsart}
\usepackage{calc,amssymb,amsmath,amsthm,amsfonts,ifthen,graphics,enumerate}
\usepackage[all]{xypic}
\newcommand{\topplace}{section}

\newcommand{\thmplace}{\topplace}
\theoremstyle{plain}
\newtheorem{theorem}{Theorem}[\thmplace]
\newtheorem{corollary}[theorem]{Corollary}

\newtheorem{proposition}[theorem]{Proposition}
\newtheorem{lemma}[theorem]{Lemma}

\newtheorem{addendum}{Addendum}[\thmplace]
\theoremstyle{definition}
\newtheorem{definition}{Definition}[\thmplace]
\newtheorem{exercise}{Exercise}[section]

\newtheorem{example}{Example}[section]

\newtheorem{remark}[theorem]{Remark}
\newtheorem{question}[exercise]{Question}

\newcommand{\ignore}[1]{\relax}
\newcommand{\ds}{\displaystyle}
\newcommand{\defem}[1]{%
  \ifmmode{\ #1\ }\else{{\mdseries\itshape\sffamily #1}}\fi\index{#1}}

\newcommand{\N}{\mathbb{N}}
\newcommand{\Tau}{\mathcal{T}}

\newcommand{\ident}{\ensuremath{\approx}}

\newcommand{\Beta}{\mathcal{B}}
\newcommand{\Assoc}{\textup{Assoc}}

\newlength{\caretsize}
\newcommand{\caret}[4]{
\xy <0mm,0mm>;<\caretsize,0mm>:
(-#1,-#2)="l";  
(0,0)="t";      
(#1,-#2)="r";   
"l";"t"**@{-};
    "r"**@{-};
"l"*!UC=(#1,#2){#3};
"r"*!UC=(#1,#2){#4};
\endxy}

\newcommand{\str}[2][10pt]{\raisebox{-#1}{$#2$}}

\newcommand{\lxzero}{\raisebox{20pt}[20pt][20pt]{
\caret{5}{10}{\caret{5}{10}{\cdot}{\cdot}}{\cdot}}}

\newcommand{\rxzero}{\raisebox{20pt}[20pt][20pt]{
\caret{5}{10}{\cdot}{\caret{5}{10}{\cdot}{\cdot}}}}

\newcommand{\lxone}{\makebox[45pt][l]{\raisebox{25pt}[25pt][25pt]{
\caret{5}{10}{\cdot}
  {\caret{5}{10}
    {\caret{5}{10}{\cdot}{\cdot}}{\cdot}}}}}

\newcommand{\rxone}{\makebox[45pt][l]{\raisebox{25pt}[25pt][25pt]{
\caret{5}{10}{\cdot}
  {\caret{5}{10}{\cdot}
    {\caret{5}{10}{\cdot}{\cdot}}}}}}

\newcommand{\lczero}{\raisebox{30pt}[25pt][25pt]{
\caret{10}{10}
  {\caret{5}{10}{\cdot}{\caret{5}{10}{\cdot}{\cdot}}}
  {\caret{5}{10}{\cdot}{\cdot}}}}

\newcommand{\rczero}{\raisebox{30pt}[25pt][25pt]{
\caret{10}{10}
  {\caret{5}{10}{\cdot}{\cdot}}
  {\caret{5}{10}{\caret{5}{10}{\cdot}{\cdot}}{\cdot}}}}

\newcommand{\lcone}{\raisebox{35pt}[35pt][30pt]{
\caret{12}{10}
  {\caret{5}{10}{\cdot}{\cdot}}
  {\caret{5}{10}
    {\caret{5}{10}
      {\caret{5}{10}{\cdot}{\cdot}}{\cdot}}{\cdot}}}}

\newcommand{\rcone}{\raisebox{35pt}[35pt][30pt]{
\caret{12}{10}
  {\caret{5}{10}{\cdot}{\cdot}}
  {\caret{5}{10}
    {\caret{5}{10}{\cdot}
      {\caret{5}{10}{\cdot}{\cdot}}}{\cdot}}}}

\newcommand{\plfunc}[1]{\xy <#1,0mm>:
(0, 0);( 4, 8)**@{-};
       ( 8,12)**@{-};
       ( 9,14)**@{-};
       (16,16)**@{-};
\endxy}

\newcommand{\bplxzero}[1]{\xy <#1,0mm>:
(0,0);( 4, 8)**@{-};
      ( 8,12)**@{-};
      (16,16)**@{-};
\endxy}

\newcommand{\plczero}{\xy <2mm,0mm>:
( 8, 8)*\bplxzero{1mm}*\frm{.};
( 0, 0);( 4, 4)**@{-};
(12,12);(16,16)**@{-};
(0,0);(16,16)**\frm{-};
(8,5)*{x_0};
\endxy}

\newcommand{\plcone}{\xy <2mm,0mm>:
(10,10)*\bplxzero{0.5mm};
( 0, 0);( 8, 8)**@{-};
(12,12);(16,16)**@{-};
( 4, 4);(12,12)**\frm{.};
( 0, 0);(16,16)**\frm{-};
(8,5)*{x_1};
\endxy}

\begin{document}


\title{ Associativity and Thompson's Group}
\author{Ross Geoghegan and Fernando Guzm\'an}

\begin{abstract}

Given a set $S$ equipped with a binary operation (we call this a
``bracket algebra") one may ask to what extent the binary operation
satisfies some of the consequences of the associative law even when it
is not actually associative?  We define a subgroup \Assoc($S$) of
Thompson's Group $F$ for each bracket algebra $S$, and we interpret
the size of \Assoc($S$) as determining the amount of associativity in
$S$ - the larger \Assoc($S$) is, the more associativity holds in $S$.
When $S$ is actually associative, \Assoc($S$) = $F$; that is the
trivial case.  In general, it turns out that only certain subgroups of
$F$ can occur as \Assoc($S$) for some $S$, and we describe those
subgroups precisely.  We then explain what happens in some familiar
examples: Lie algebras with the Lie bracket as binary operation,
groups with the commutator bracket as binary operation, the Cayley
numbers with their usual multiplication, as well as some less familiar
examples.  In the case of a group $G$, with the commutator bracket as
binary operation, it is better to think of the ``virtual size of $G$",
determined by all the groups \Assoc($H$) such that $H$ is a subgroup
of finite index in $G$.  This gives a way of partitioning groups into
``small", ``intermediate" and ``large" - a partition suggestive of,
but different from, traditional measures of a group's size such as
growth, isoperimetric inequality and ``amenable vs. non-amenable"

\end{abstract}

\maketitle

\section{Introduction}

\subsection{Thompson's Group}

We begin by recalling the ``pairs of binary trees" definition of
Thompson's group $F$ as we will need precise terminology.

By a \defem{binary tree} we will mean either the trivial tree
consisting of a single vertex and no edges, or a finite tree having
one vertex of order 2 (called the \defem{root}) and all other vertices
having order either 1 (they are the \defem{leaves}) or 3 (they are the
\defem{interior vertices}).  The only vertex in the trivial tree is
considered to be both the root and a leaf.  Included in the structure
of what we call a binary tree is a labeling of the edges as follows:
each edge has a canonical orientation away from the root, and in
non-trivial trees exactly two edges point away from each non-leaf
vertex; one of those edges is labeled $0$ and the other $1$.  This is
most intuitively seen when the binary tree is embedded in the
Euclidean plane as in the pictures below, with all oriented edges
pointing downward, $0$ to the left and $1$ to the right.  Indeed,
every binary tree (in our sense) can be represented in this way.  We
will not distinguish between two binary trees when there exists a
label-preserving isomorphism between them.

The leaves of a binary tree have a canonical ordering as follows:
there is a unique geodesic edge path from the root to each leaf
specified by a word in the alphabet $\{0,1\}$: $0$ means follow the
left edge, $1$ means follow the right edge.  For example, in the
figure below, the leaf marked $b$ is specified by $010$ and the leaf
marked $d$ is given by $10$; The canonical ordering of the leaves is
then the corresponding lexicographic ordering of those geodesics.  By
the $i$th leaf we will mean the $i$th term in this ordering; if there
are $n$ leaves they are numbered $1,\dots,n$.

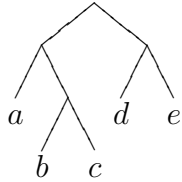
\[
\raisebox{35pt}[20pt][30pt]
{
\caret{10}{8}
  {\caret{5}{10}{\str a}{\caret{5}{10}{\str b}{\str c}}}
  {\caret{5}{10}{\str d}
    {\str e}}}
\]

A \defem{caret} is a binary tree with exactly two edges.  If $p$ is a
binary tree having $n$ leaves, the $i$-th \defem{elementary expansion}
of $p$ is the binary tree $\beta^i(p)$ obtained from the disjoint
union of $p$ and a single caret by identifying the root of the caret
with the $i$th leaf of $p$.  This makes sense when $1\leq i \leq n$;
it is convenient to define $\beta^i(p) = p$ when $i > n$.

Let $\Tau_n$ denote the set of (isomorphism classes of) binary trees
which have $n$ leaves.  We write $ \Tau^{(2)} =
\coprod_{n=1}^{\infty}{\Tau_n \times \Tau_n}$ and we define $b^i:
\Tau^{(2)} \to \Tau^{(2)}$ to agree with $\beta^i \times \beta^i$ on
$\Tau_n \times \Tau_n$.  The relations $(p,q) \sim b^i(p,q)$, one such
relation for each $i$, together generate an equivalence relation on
$\Tau^{(2)}$.  The set of equivalence classes is denoted by $F$.  We
write $\langle p,q\rangle$ for the equivalence class of $(p,q)$

To define a multiplication on $F$ we note first that any two binary
trees have a common expansion, where an \defem{expansion} of $p$ is a
binary tree of the form $\beta(p)= \beta^{i_k}\dots \beta^{i_1}(p)$.
We define a \defem{simultaneous expansion} of the ordered pair $(p,q)$
to be an ordered pair of the form $(\beta(p), \beta(q))$.  The
\defem{product} $\langle p,q\rangle.\langle r,s\rangle$ is then
defined to be $\langle p',s'\rangle$, where $(p',q')$ is an expansion
of $(p,q)$, $(r',s')$ is an expansion of $(r,s)$, and $q' = r'$.  This
is well-defined and associative.  With respect to this multiplication
the element $\langle p,p\rangle$ (where $p$ is any binary tree - this
is independent of $p$) is a two-sided identity, and $\langle
q,p\rangle$ is a multiplicative inverse for $\langle p,q\rangle$.
Thus we have a group - Thompson's group $F$.

The group $F$ has \defem{standard generators} $x_0$ and $x_1$
represented by the following two pairs of binary trees:

\[
\setlength{\caretsize}{0.5mm}
\left( \lxzero \ \ , \rxzero\ \right)
\]

\[
\setlength{\caretsize}{0.5mm}
\left( \lxone ,\rxone \right)
\]

For details, and for information about some of the many remarkable
properties of this group, see \cite{BG} or the expository article
\cite{CFP}.


\subsection{Bracket Algebras}

By a \defem{bracket algebra} we will mean a set $S$ together with a
binary operation $\alpha:S \to S$.  A binary tree having $n$ leaves
determines a rule for associating ordered $n$-tuples of members of $S$
using $\alpha$, i.e. determines an $n$-ary operation on $S$; see
Section~\ref{section two} for details.  So an ordered pair $(p,q)$ of
binary trees having the same number of leaves, can be interpreted as
encoding a ``law" in $S$ saying that the method of association defined
by $p$ always gives the same result in $S$ as the method of
association defined by $q$.  In that case we can say that $S$
``satisfies" $(p,q)$.  Similarly, an element $\langle p,q\rangle$ of
$F$ can be interpreted as defining a ``stable law" in $S$; i.e. $S$
satisfies some simultaneous expansion of $(p,q)$.  Note that if $S$
satisfies $(p,q)$ then it satisfies every simultaneous expansion of
$(p,q)$.

\begin{proposition}\label{prop:subgroupF} 
  The elements $\langle p,q\rangle$ of $F$ which define stable laws in
  the bracket algebra $S$ form a subgroup of $F$.
\end{proposition}

The proof is in Section 2.  We call the subgroup in
Proposition~\ref{prop:subgroupF} \Assoc($S$), the \defem{group of
  stable associativities} of $S$.


\subsection{Characterization of the groups \Assoc($S$)}

A natural first question is: which subgroups of $F$ can occur as
\Assoc($S$) for some bracket algebra $S$?  To answer this we need some
further vocabulary concerning binary trees.

The \defem{right shift} of $p$ is the binary tree $\sigma_1(p)$
obtained from the disjoint union of $p$ and a single caret by
identifying the root of $p$ with the second (i.e. right) leaf of the
caret.  The \defem{left shift} of $p$, $\sigma_0(p)$, is defined
similarly.

\[ \sigma_0(p)= \raisebox{15pt}{\caret{5}{10}{\str[8pt]p}{\cdot}}\quad 
\sigma_1(p)= \raisebox{15pt}{\caret{5}{10}{\cdot}{\str[8pt]p}} \]

There is an involution on the set of binary trees taking a binary tree
$p$ to its \defem{reflection} $\rho(p)$.  In terms of planar pictures,
$\rho(p)$ is the mirror image of $p$ in the $y$-axis.

The two shifts induce endomorphisms $s_i:F \to F$ ($i=$ 0 or 1) taking
$\langle p,q\rangle$ to $\langle \sigma_i(p),\sigma_i(q)\rangle$.
These are the \defem{left and right shift endomorphisms}.  The
\defem{reflection} automorphism $R:F \to F$ is defined by $R(\langle
p,q\rangle) = \langle \rho(p), \rho(q)\rangle$.  Note that $s_0 =
R.s_{1}. R$.  The right shift $s_1$ is, of course, well-known, and is
defined by the formula $s_1(x_0) = x_1$ and $s_1(x_1) = x_2$ where
$x_2 := x_1^{x_0}$.  The formula for the left shift in terms of the
standard generators is less pleasant: $s_0(x_0) =
{(x_{0}x_{1}^{-1})}^{x_0^{-1}}$ and $s_0(x_1) =
{(x_{1}x_{2}^{-1})}^{(x_{0}x_{1})^{-1}}$.

The shift endomorphisms can be much better understood in terms of the
``dyadic piecewise linear" model of $F$. This is a well-known faithful
representation of $F$ in the group of increasing self-homeomorphisms
of the closed unit interval $I$.  Given $\langle p,q\rangle \in F$,
two copies of $I$ are to be dyadically subdivided according to the
instructions of the trees $p$ and $q$, and then $\langle p,q\rangle$
is identified with the dyadic piecewise linear increasing
homeomorphism of $I$ which maps each segment of the $p$-subdivision
affinely onto the corresponding segment of the $q$-subdivision.  In
these terms, the two shift endomorphisms can be easily understood from
the following picture:

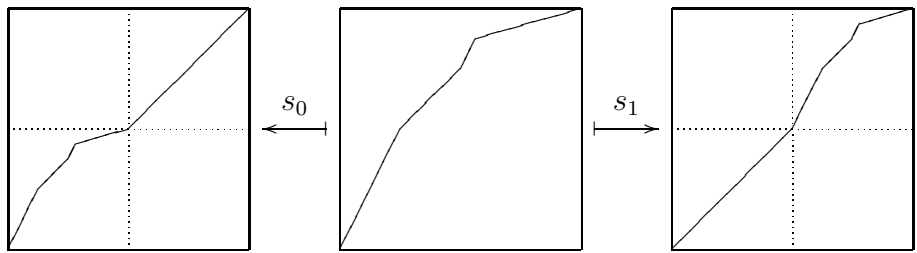
\[\xymatrix{
\xy <1mm,0mm>:
(0,0);(32,32)**\frm{-};
( 8, 8)*\plfunc{1mm}*\frm{.};
(16,16);(32,32)**@{-}**\frm{.};
\endxy \ar@{<-|}[r]^{\mbox{\normalsize $s_0$}}
&\xy <2mm,0mm>:
(8,8)*\plfunc{2mm}*\frm{-}
\endxy \ar@{|->}[r]^{\mbox{\normalsize $s_1$}}
&\xy <1mm,0mm>:
(0,0);(32,32)**\frm{-};
(24,24)*\plfunc{1mm}*\frm{.};
( 0, 0);(16,16)**@{-}**\frm{.}
\endxy} \]

Now we can answer the question about what subgroups of $F$ occur as
\Assoc($S$):

\begin{theorem}\label{theorem:characterization} 
  A subgroup $H$ of $F$ is \Assoc($S$) for some bracket algebra $S$ if
  and only if both shift endomorphisms, $s_0$ and $s_1$, map $H$ into
  itself.
\end{theorem} 

The proof is at the end of Section 2.  Examples~\ref{exam:not normal}
and~\ref{exam:small not normal} illustrate how a non-normal subgroup
of $F$ can arise as \Assoc($S$).  There are infinitely many such
subgroups.

A companion question is: which normal subgroups of $F$ occur as
\Assoc($S$) for some bracket algebra $S$?  Recall that $F/F'$ is
isomorphic to $\mathbb{Z} \times \mathbb{Z}$ where the standard
generators $x_0$ and $x_1$ map respectively to $(1,0)$ and
$(0,1)$. Since the non-trivial normal subgroups of $F$ are the
subgroups containing $F'$, they are in
bijective correspondence with the subgroups of $\mathbb{Z} \times
\mathbb{Z}$.  Both shifts preserve $F'$ and so induce endomorphisms
$\bar{s}_i$ of $\mathbb{Z} \times \mathbb{Z}$; the formulas are:
$\bar{s_0}(m,n) = (m,-m)$ and $\bar{s_1}(m,n) = (0,m+n)$.  Thus we
have:

\begin{corollary}\label{cor:normal characterization} 
  A normal subgroup of $F$ is \Assoc($S$) for some $S$ if and only if
  it is the preimage in $F$ of a subgroup of $\mathbb{Z} \times
  \mathbb{Z}$ generated by $(m,-m)$ and $(0,n)$ for some integers $m$
  and $n$.
\end{corollary}

In Proposition~\ref{prop:normal} we give a bracket algebra criterion
for \Assoc($S$) to be normal.


\subsection{Strongly Regular Laws}

Given a bracket algebra $S$, an element $(p,q)$ of $\ds\Tau^{(2)}$ may
or may not hold as a law in $S$. For example, the associative law in
$S$
\[ \alpha(\alpha(x_1,x_2),x_3) \ident \alpha(x_1,\alpha(x_2,x_3)) \]
corresponds to the pair of trees 
{\setlength{\caretsize}{0.5mm}
\[ \left( \lxzero\ , \rxzero \ \ \right) \]
We write  $S\models p\ident q$ when that law holds in $S$. 

Such laws are of a special type: each variable occurs exactly once on
each side, and variables occur in the same order on both sides.  We
will call such laws (represented by pairs $(p,q)\in\Tau^{(2)}$)
\defem{strongly regular laws}.  Similarly, we say that the law
$p\ident q$ \defem{eventually holds} in $S$, and write $S\models_e
p\ident q$, if there is a simultaneous expansion $(p',q')$ of $(p,q)$
such that $S\models p'\ident q'$.  In these terms we obviously have:

\begin{proposition}\label{prop:large is lawless} 
  \Assoc($S$) is trivial if and only if no strongly regular law
  eventually holds in $S$.
\end{proposition}


\subsection{Example: Lie algebras}

For Lie algebras with the Lie bracket as the binary operation we have:

\begin{theorem}\label{thm:fd lie algebras statement} 
  Let $L$ be a finite-dimensional complex Lie algebra. \Assoc($L$) is
  either trivial or is $F$.  It is $F$ if and only if $L$ is a
  solvable Lie algebra.
\end{theorem} 

The proof is in Section 4.  

Kac-Moody algebras are infinite-dimensional, but the proof of
Theorem~\ref{thm:fd lie algebras statement} will show:

\begin{addendum}\label{Kac-Moody addendum} 
  If $L$ is a Kac-Moody algebra then \Assoc($L$) is trivial.
\end{addendum}


\subsection{Example: Groups}

A group $G$ with the commutator bracket $[x,y] = xyx^{-1}y^{-1}$ as
the binary operation is an important example of a bracket algebra.
When we discuss $G$ as a bracket algebra this binary operation is
always understood.

\begin{proposition}\label{prop:groups statement} 
  Let $G$ be a group. \Assoc($G$) = $F$ if and only if $G$ is a
  solvable group.
\end{proposition}

\begin{proposition}\label{prop:A5} 
  \Assoc($A_5$) is trivial, where $A_5$ denotes the finite alternating
  group on five letters.
\end{proposition}

Proposition~\ref{prop:groups statement} is proved in Section 3, and
Proposition~\ref{prop:A5} is proved in Section 5.
Proposition~\ref{prop:A5} is of interest because $A_5$ has a subgroup
of finite index, the trivial group $\{1\}$, such that \Assoc($\{1\}$)
is $F$.  Such instability suggests that in order to use \Assoc($G$) as
a good measure of the size of a group $G$ we should consider things
virtually.

We define $G$ to be \defem{ large} if for every subgroup $H$ of finite
index \Assoc($H$) is trivial.  Equivalently, by
Proposition~\ref{prop:large is lawless}, $G$ is large if and only if
no strongly regular law holds eventually in any subgoup $H$ of finite
index in $G$.

We define $G$ to be \defem{small} if it has a subgroup $H$ of finite
index such that Assoc($H$) is a non-trivial normal subgroup of $F$.
(Recall that a subgroup of $F$ is normal if and only if it is either
trivial or contains the commutator subgroup $F'$.)

If $G$ satisfies the third possibility, neither large nor small, we
say that $G$ is \defem{of intermediate size}.  We do not know the
answer to the following:

\begin{question} 
  Does there exist a group of intermediate size?
\end{question} 

Indeed, we do not know a group $G$ for which \Assoc($G$) is neither
$F$ nor the trivial group.

Let $F_k$ denote the free group on $k$ generators $\{X_1, \dots,X_k
\}$ and let $w(X_1, \dots,X_k)$ be a non-empty word in those
generators and their inverses.  This word $w$ defines a \defem{law} in
a group $G$ if the statement $w(g_1,\dots,g_k) = 1$ is true for all
choices of $g_i \in G$.  A group \defem{satisfies no law stably} if
each subgroup of finite index satisfies no law.  Every strongly
regular law in the commutator algebra of $G$ can be recast as a law in
this sense, so by Proposition~\ref{prop:large is lawless} every group
which satisfies no law stably is large.  This and some related
observations are summarized in:

\begin{proposition}\label{prop:general stuff} 
  Groups satisfying no law stably are large.  If a group $G$ has a
  quotient containing a large subgroup then $G$ is also large.
  Virtually solvable groups are small.  
\end{proposition}

This is proved in Section 3.  

Many interesting groups are known to satisfy no law stably, for
example non-abelian free groups, non-elementary hyperbolic groups, and
Thompson's Group $F$, and are therefore large. (Proof of this for $F$:
that $F$ itself satsifies no law was proved in \cite{BS}; every
subgroup of finite index in $F$ contains a normal subgroup of finite
index, which in turn contains the commutator subgroup, and hence
contains a copy of $F$.)  Moreover, by Proposition~\ref{prop:general
  stuff} the Tits Alternative separates the finitely generated linear
groups neatly into the large and the small (i.e. every such group
either contains a free non-abelian subgroup or is virtually solvable).

The relationship between ``stably large" and ``non-amenable" is not so
clean.  Among the large groups are the non-abelian free groups (which
are non-amenable) and Thompson's group $F$ (which is conjectured to be
non-amenable).  However, it is not true that ``non-amenable" is
equivalent to ``stably large".  An example is the ``first Grigorchuk
group" $\Gamma$ (\cite{Gr} or \cite{BGS}), which is finitely
generated, has subexponential growth (and is therefore amenable) but
is large since it satisfies no law stably (\cite{A} or \cite{S}).  Now
this group $\Gamma$ is not elementary amenable, but even in the
elementary amenable case there are distinctions to be made.  We ask:

\begin{question}\label{elementary amenable} 
  Is there a finitely generated large elementary amenable group?
\end{question}

The words ``finitely generated" are included in this question because
of:

\begin{proposition}\label{proposition:Sinfinity} 
  Let $S_n$ denote the group of permutations of $\{1, \dots, n\}$ and
  let $S_{\infty}$ denote the union of the groups $S_n$.  Then
  $S_{\infty}$ is elementary amenable and is large.
\end{proposition}

This is proved in Section 5.


\subsection{The Five Variable Law}

The commutator of the two standard generators $[x_0, x_1] =
x_{0}x_{1}x_{0}^{-1}x_{1}^{-1}$ is most efficiently represented by the
following pair of binary trees.

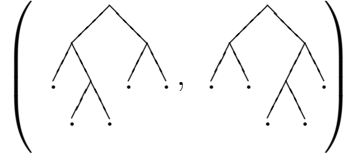
\[ \left( \lczero\ , \rczero \right) \]

This pair of trees defines a strongly regular law which we call the
Five Variable Law.

\begin{theorem}\label{thm:eventual FVL} 
  The commutator $[x_0, x_1]$ lies in \Assoc($S$) if and only if
  \Assoc($S$) is a non-trivial normal subgroup of $F$.  In other
  words, \Assoc($S$) is a non-trivial normal subgroup of $F$ if and
  only if the Five Variable Law holds eventually in $S$.
\end{theorem}

This is proved in Section 6.

An elegant application is that the Five Variable Law holds eventually
in a finite-dimensional complex Lie algebra $L$ if and only if $L$ is
solvable.  Perhaps this is well-known?

We call the bracket algebra $(S, \alpha)$ \defem{simply perfect} if
$\alpha$ is surjective.

\begin{remark}\label{rem:simply perfect} 
  It is clear that if $S$ is simply perfect and if some simultaneous
  expansion of the Five Variable Law holds in $S$ then the Five
  Variable Law itself holds in $S$. Indeed, this holds for any
  strongly regular law.  This immediately implies the next Corollary.
\end{remark}

\begin{corollary}\label{cor:FVL on the nose} 
  Let $S$ be a simply perfect bracket algebra.  \Assoc($S$) is a
  non-trivial normal subgroup of $F$ if and only if the Five Variable
  Law holds (``on the nose") in $S$.  
\end{corollary}


\subsection{Example: The Cayley numbers; bracket algebras having
  identity elements} 

A classical example of a non-associative bracket algebra is the
algebra of Cayley numbers, also known as the octonion algebra; its
underlying real vector space is eight-dimensional.  Since it possesses
a two-sided identity element $1$, it is covered by the following
general theorem:

\begin{theorem}\label{thm:identity element} 
  Let $S$ be a bracket algebra which has a two-sided identity element
  $1$.  If $S$ is non-associative then it satisfies no non-trivial
  strongly regular identity.  Thus \Assoc($S$) is either $F$ or is
  trivial depending on whether $S$ is associative or not.
\end{theorem}

In particular, this theorem (which is proved in Section 7) applies to
all loops.  (A \defem{loop} is a bracket algebra with two- sided
identity element such that each element has a unique right inverse and
a unique left inverse.)

By contrast, we exhibit a four-element non-associative bracket algebra
$S$ having a right identity element but no left identity element such
that \Assoc($S$) is a non-normal subgroup of $F$; see
Example~\ref{exam:small not normal}.


\subsection{Acknowledgments}

An early version of some of these ideas was presented by the
first-named author in the symposium ``Thompson's Group at Forty Years"
at the American Institute of Mathematics in January 2004.  He thanks
that institution for their support.  He also takes this opportunity to
warmly thank his colleagues for their dedication of this volume, and
for attending the accompanying special session held in his honor at
the October 2004 meeting of the American Mathematical Society at
Vanderbilt University.  And he expresses appreciation to Michael
Mihalik and Mark Sapir who organized that session.

\section{Bracket algebras and Thompson's Group}\label{section two}

Let $\Tau$ denote the set of all binary trees, made into a bracket
algebra as follows.  The trivial binary tree will be denoted by
$\cdot$, and the binary tree having $p$ as the left subtree and $q$ as
the right subtree,
\[  \raisebox{15pt}{\caret{5}{10}{\str[8pt]p}{\str[8pt]q}} \]
\noindent
will be denoted by $[p,q]$.  This is the free bracket algebra
generated by $\cdot$, the trivial binary tree.  Recall that we
extended the function $\beta^i: \Tau_n \to \Tau_{n+1}$ to a function
$\Tau \to \Tau$ by setting $\beta^i(p)=p$ for $p\in\Tau_n$ when $i>n$.
Let $\Beta$ denote the monoid generated by $\{\beta^i|i\in\N\}$ with
composition as multiplication (the \defem{monoid of expansions}).
Thus we have an action of $\Beta$ on the set $\Tau$.  Observe that the
elementary expansions $\beta^i$, for $i\in\N$ satisfy: \[
\beta^i\beta^j = \beta^{j+1}\beta^i, \text{ when } i < j. \] (Indeed,
it can be shown that $\Beta$ is isomorphic to the negative monoid of
Thompson's group $F$.)

Now, let $\alpha:S\times S\to S$ be a bracket algebra.  (We will
sometimes denote the binary operation in such an algebra by $[,]$
rather than by $\alpha$, and we will often omit explicit reference to
the operation altogether.)  We define a bracket algebra $X(S)$, a
quotient $\Beta$-set of the $\Beta$-set $\Tau$, which makes precise
the notion of ``expansion" in the context of any bracket algebra.

First we define $X_n(S)$ and a bracket on $\coprod_{n=1}^{\infty}
X_n(S)$ recursively by:

\begin{itemize}
\item The identity map, $I_S:S\to S$, is the only member of $X_1(S)$
\item if $f\in X_n(S)$, $g\in X_m(S)$, and
    $\overline{s}=(s_1,\dots,s_n,s_{n+1},\dots,s_{n+m})$, then
    \[\begin{array}{rccl}
    $[$f$,$g$]$:&S^{n+m}&\to &S \\
                &\overline{s}&\mapsto
                    &\alpha(f(s_1,\dots,s_n),g(s_{n+1},\dots,s_{n+m})) \\
    \end{array} \]
    belongs to $X_{n+m}(S)$.
\end{itemize}

The promised $X(S)$ is $\coprod_{n=1}^{\infty} X_n(S)$ with this
bracket.  It is the bracket algebra of ``$\alpha$-operations" on $S$,
generated by one element, namely, $I_S$.  Hence, it is a quotient of
$\Tau$, with quotient epimorphism $\widehat{\ }:\Tau\to X(S)$ given
recursively by
\begin{itemize}
\item $\hat{\cdot}=I_S$ where $\cdot$ denotes the trivial binary tree, and
\item $\widehat{[p,q]}(s_1,\dots,s_n,s_{n+1},\dots,s_{n+m})=
  \alpha(\widehat{p}(s_1,\dots,s_n),\widehat{q}(s_{n+1},\dots,s_{n+m}))$
\end{itemize}

\begin{proposition}\label{prop:order}
  The action of $\Beta$ on $\Tau$ descends via the epimorphism
  $\widehat{\ }:\Tau\to X(S)$ to an action of $\Beta$ on $X(S)$, so
  that $\widehat{\ }$ is also a morphism of $\Beta$-sets.
\end{proposition}

\begin{proof}
  Clearly the function $\widehat{\ }$ is a surjective morphism of
  bracket algebras.  If $p,q\in\Tau$ are such that
  $\widehat{p}=\widehat{q}$ then $p$ and $q$ have the same number $n$
  of leaves.  For $i\leq n$
  \[ \begin{array}{rcl}
    \widehat{\beta^i(p)}(s_1,\dots,s_i,s_{i+1},\dots,s_n)
    &=&\widehat{p}(s_1,\dots,\alpha(s_i,s_{i+1}),\dots,s_n) \\
    &=&\widehat{q}(s_1,\dots,\alpha(s_i,s_{i+1}),\dots,s_n) \\
    &=&\widehat{\beta^i(q)}(s_1,\dots,s_i,s_{i+1},\dots,s_n) \\
  \end{array}
  \]
  so, $\widehat{\beta^i(p)}=\widehat{\beta^i(q)}$. \\
\end{proof}  

For $\beta\in\Beta$ we call $\beta(f)$ an \defem{expansion} of $f \in
X(S)$.

The action of $\Beta$ on $\Tau$ induces a diagonal action of $\Beta$
on $\Tau^{(2)}$ and the set of orbits, $\Tau^{(2)}/ \Beta$, is
precisely Thompson's Group $F$.

We write $X_n^{(2)}(S)=X_n(S)\times X_n(S)$ and $\ds
X^{(2)}=\coprod_{n=1}^\infty X_n^{(2)}$.  We call
$(\beta(f),\beta(g))$ a \defem{simultaneous expansion} of the pair
$(f,g)$.  Again using the diagonal action of $\Beta$, we consider the
set of orbits $X^{(2)}/ \Beta$.

Clearly we have: 

\begin{proposition}\label{prop:assoc is kernel}
  \begin{enumerate}
  \item\label{prop:assoc is kernel part:compatible}
    The function
    \[ \begin{array}{rrcl}
      \widehat{\ }:&\Tau^{(2)}&\to &X^{(2)} \\
                   &(p,q) &\mapsto &(\widehat{p},\widehat{q}) \\
    \end{array}
    \]
    is a surjective $\Beta$-map.  
  \item
    The pre-image under the induced function  
    \[ \begin{array}{rrcl}
      \widehat{\ }:&\Tau^{(2)}/\Beta &\to &X^{(2)}/\Beta \\
       &\langle p,q\rangle &\mapsto &\langle\widehat{p},\widehat{q}\rangle \\
    \end{array}
    \]
    of the trival element of $X^{(2)}/\Beta$ (i.e. of the image of $1
    \in F = \Tau^{(2)}/\Beta$) is the set \Assoc($S$).
  \end{enumerate}
\end{proposition}

When we consider strongly regular identities which eventually hold in
$S$, we look at pairs $(p,q)\in \Tau_n\times\Tau_n$, such that
$S\models_e p\ident q$, i.e. such that $\widehat{p'}=\widehat{q'}$
where $(p',q')=(\beta(p), \beta(q))$ for some simultaneous expansion
$\beta \times \beta$.  Of course, this is the same as looking at
classes $\langle p,q\rangle \in \Tau^{(2)}/\Beta (=F)$.

Consideration of binary trees leads immediately to the observation
that $\Beta$ has the ``common left multiples property'', i.e. given
$\beta_1,\beta_2\in\Beta$ there are $\beta_3,\beta_4\in\Beta$ such
that $\beta_3\beta_1=\beta_4\beta_2$.

Proposition~\ref{prop:subgroupF}  can be restated as:

\begin{proposition}\label{prop:assoc}
  Let $S$ be a bracket algebra.  The set of stable associativities
  $\Assoc($S$)= \{\langle p,q\rangle \in F \ | S\models_e p\ident q
  \}$ is a subgroup of Thompson's Group $F$.
\end{proposition}

\begin{proof}
  Clearly \Assoc($S$) contains the identity element of $F$ and is
  closed under inverses.  For closure under composition, let $\langle
  p,q\rangle ,\langle s,r\rangle \in F$.  WLOG we may assume $q=s$, so
  we have $S\models_e p\ident q$ and $S\models_e q\ident r$.  Let
  $\beta_1,\beta_2$ be expansions such that
  $\beta_1(\widehat{p})=\beta_1(\widehat{q})$ and
  $\beta_2(\widehat{q})=\beta_2(\widehat{r})$.  By ``common left
  multiples" there are expansions $\beta_3,\beta_4$ such that
  $\beta_3\beta_1=\beta_4\beta_2$.  Let
  $\beta=\beta_3\beta_1=\beta_4\beta_2$. Then,
  $\beta(\widehat{p})=\beta_3\beta_1(\widehat{p})=\beta_3\beta_1(\widehat{q})
  =\beta_4\beta_2(\widehat{q})=\beta_4\beta_2(\widehat{r})=\beta(\widehat{r})$.
  So, \mbox{$S\models_e p\ident r$}.
\end{proof}

The next proposition gives us necessary and sufficient conditions for
\Assoc($S$) to be a normal subgroup of $F$.  First a definition.

\begin{definition}\label{def:normal}
  Let $S$ be a bracket algebra.  We say that $S$
  is \defem{normal} provided that: whenever $\beta_1,\beta_2\in\Beta$
  and $g\in X(S)$ are such that $g\neq\beta_1(g)=\beta_2(g)$, then for
  any $f\in X(S)$ there is $\beta\in\Beta$ such that
  $\beta\beta_1(f)=\beta\beta_2(f)$. 
\end{definition}

\begin{proposition}\label{prop:normal}
  Let $S$ be a bracket algebra.  $\Assoc($S$)$ is a normal subgroup of
  $F$ if and only if $S$ is a normal bracket algebra.
\end{proposition}

\begin{proof}
  Assume first that the bracket algebra $S$ is normal.  Let $\langle
  p,q \rangle \in$ \Assoc($S$).  WLOG we may assume that $\hat{p} =
  \hat{q}$.  For any $\langle r,s \rangle$ we have

  \[ \begin{array}{rcl}
    \langle r,s \rangle \cdot \langle p, q \rangle \cdot \langle s,r
    \rangle &=& \langle{\beta_{1}(r)}, \beta_{1}(s) \rangle \cdot
    \langle \beta_{2}(p), \beta_{2}(q) \rangle \cdot \langle
    \beta_{3}(s), \beta_{3}(r) \rangle\\ &=&\langle \beta_{1}(r),
    \beta_{3}(r) \rangle\\
  \end{array}
  \]

  where the $\beta_i$ are chosen so that $\beta_1(s) = \beta_{2}(p)$
  and $\beta_{2}(q) = \beta_{3}(s)$.

  So $\widehat{\beta_1(s)}= \widehat{\beta_2(p)}=
  \widehat{\beta_2(q)}= \widehat{\beta_3(s)}$.  Thus, either $s\neq
  \beta_1(s) = \beta_{3}(s)$ in which case, by normality, there exists
  $\beta$ such that $\widehat{\beta \beta_1(r)}= \widehat{\beta
    \beta_3(r)}$, implying that our conjugate of $\langle p,q \rangle$
  also lies in \Assoc($S$); or $s=\beta_{1}(s) = \beta_{3}(s)$ in
  which case $\beta_1(r)=\beta_3(r)=r$, which would mean that $\langle
  p,q \rangle$ is the trivial element of $F$.

  Conversely, if \Assoc($S$) is normal then the multiplication on
  $F=\widehat{\Tau^{(2)}}/\Beta$ descends to $X^{(2)}/\Beta$ and the
  condition for normality of $S$ follows easily. \\
\end{proof}

We now describe an example where the group \Assoc($S$) is not normal
in $F$.  For this, we need further notation concerning binary trees.

If $p$ is a binary tree, with each vertex of $p$ is associated a
unique word $w$ in the alphabet $\{0,1\}$.  This was explained for
leaves in the Introduction but applies equally to all vertices.  The
leaves of $p$ are numbered $1$ through $n$ from left to right, and the
word of the $i$-th leaf is denoted by $l_i$ or $l_i(p)$.  We say that
$p$ has a \defem{free caret} at the interior vertex $w$ if both $w0$
and $w1$ denote leaves of $p$.  The \defem{depth} of the $i$-th leaf
is the length of the word that denotes it, $d_i=|l_i|$.  If $w$
denotes a vertex of $p$, the subtree of $p$ at that vertex is denoted
by $p_w$.

Here is the promised example:

\begin{example}\label{exam:not normal}
  Let

  \[ p=\lxone \quad q=\rxone \]
  \vspace{10pt}

\noindent
  Let $S$ be a free algebra in the variety of bracket algebras that
  satisfy the identity: $p\ident q$.
  Consider the expansions  $\beta^4$, $\beta^2$ and the 
  function $g(X,Y,Z,W)=\widehat{p}(X,Y,Z,W)=\widehat{q}(X,Y,Z,W)$.
  Notice that $\beta^4(p)=\beta^2(q)$, so $\beta^4(g)=\beta^2(g)$. 
  Let 

\[ \setlength{\caretsize}{0.5mm}
r= \makebox[65pt][l]{\raisebox{20pt}[20pt][20pt]{
\caret{8}{8}{\caret{5}{10}{\cdot}{\cdot}}{\caret{5}{10}{\cdot}{\cdot}}}}
\hspace{30pt}
r_1=\beta^4(r)= \makebox[65pt][l]{\raisebox{20pt}[20pt][20pt]{
\caret{8}{8}
  {\caret{5}{10}{\cdot}{\cdot}}
  {\caret{5}{10}{\cdot}
    {\caret{5}{10}{\cdot}{\cdot}}}}}
\hspace{30pt}
r_2=\beta^2(r)= \makebox[65pt][l]{\raisebox{20pt}[20pt][20pt]{
\caret{8}{8}
  {\caret{5}{10}{\cdot}{\caret{5}{10}{\cdot}{\cdot}}}
  {\caret{5}{10}{\cdot}
    {\cdot}}}}
\]

\noindent
  We claim that $S$ does not satisfy any simultaneous expansion of
  $r_1\ident r_2$.  Let $\beta\in\Beta$ and let $r'_1=\beta(r_1)$, and
  $r'_2=\beta(r_2)$.   There are $t_1,\dots,t_5\in\Tau$ such that 

{\normalsize
\[ \setlength{\caretsize}{0.7mm}
r'_1=r_1(t_1,\dots,t_5)=
\makebox[65pt][l]{\raisebox{20pt}[20pt][20pt]{
\caret{10}{8}
  {\caret{5}{10}{\str{t_1}}{\str{t_2}}}
  {\caret{5}{10}{\str{t_3}}
    {\caret{5}{10}{\str{t_4}}{\str{t_5}}}}}}
\hspace{30pt}
r'_2=r_2(t_1,\dots,t_5)=
\makebox[65pt][l]{\raisebox{20pt}[20pt][20pt]{
\caret{10}{8}
  {\caret{5}{10}{\str{t_1}}{\caret{5}{10}{\str{t_2}}{\str{t_3}}}}
  {\caret{5}{10}{\str{t_4}}
    {\str{t_5}}}}}
\]}
\vspace{20pt}

\noindent
  If it were true that $S\models r'_1\ident r'_2$ then it would be
  possible to rewrite $r'_1$ to $r'_2$ by a finite number of
  applications of $p\ident q$ on different subtrees.  (Strictly, this
  informal statement is a consequence of the Completeness Theorem of
  Equational Logic.)  If $t$ is any one of the intermediate trees,
  then applying $p\ident q$ at the root of $t$ does not move any
  leaves from the left to the right subtree of $t$; and applying
  $p\ident q$ at an interior vertex of $t$ does not either.  But any
  leaf of $r'_1$ under $t_3$ is in the left subtree of $r'_1$, and any
  leaf of $r'_2$ under $t_3$ is in the right subtree of $r'_2$.
\end{example}

The pair $(p,q)$ used in this Example has the property that $\langle
p,q \rangle = x_1 \in F$.  It will follow from the proof of
Theorem~\ref{theorem:characterization}, below, that the group
\Assoc($S$) in Example~\ref{exam:not normal} is the smallest subgroup
of $F$ containing $x_1$ which is invariant under both shifts.  Viewed
in the ``dyadic piecewise linear" model this subgroup is the
stabilizer of $\{\frac{1}{2^n}|n\in {\mathbb{N}}\}$.  It is abstractly
isomorphic to a countably infinite weak direct product of copies of
$F$.

The proof of Theorem~\ref{theorem:characterization} can now be given.
The ``only if" part is clear: i.e. \Assoc($S$) is invariant under both
shifts.  The``if" part follows from the following more precise Lemma:

\begin{lemma}\label{sigmainvariant}  
  Let $K:= \{\langle p_m,q_m\rangle\mid m\in {\mathcal M}\}$ be a
  subset of $F$ and let $S$ be a free algebra in the variety defined
  by $\bar K := \{(p_m,q_m)\mid m \in {\mathcal M\}}$; i.e. the
  defining set for $\bar K$ contains one representative of each member
  of $K$.  Then {\rm Assoc}$(S)$ is the smallest subgroup $H$ of $F$
  containing $K$ and invariant under both shifts.
\end{lemma}

\begin{proof}  
  Let $\langle r,s\rangle \in \Assoc(S)$.  WLOG we may assume that
  $S\models r\ident s$.  Then we can pass from $r$ to $s$ in finitely
  many steps $r = t^0,t^1,\cdots ,t^n=s$ so that each $(t^i,t^{i+1})$
  consists of a pair of trees which are identical except that at a
  certain vertex $v$ and for some $m\in{\mathcal M}$, $t^i_v=p'_m$ and
  $t^{i+1}_v=q'_m$ (or viceversa) for a simultaneous expansion
  $(p'_m,q'_m)$ of $(p_m,q_m)\in \bar{K}$.  Let $\epsilon_1\cdots
  \epsilon_k$ be the word in the alphabet $\{0,1\}$ which labels
  $v$. For any tree $p$ define $\sigma_v(p) :=
  \sigma_{\epsilon_1}\sigma_{\epsilon_2} \cdots
  \sigma_{\epsilon_k}(p)$.  Then $\langle t^i,t^{i+1}\rangle = \langle
  \sigma_v(p'_m),\sigma_v(q'_m)\rangle$ (or $\langle
  \sigma_v(q'_m),\sigma_v(p'_m)\rangle$ in the ``vice versa'' case) in
  $F$.  Thus $\langle r,s\rangle = \langle t^0,t^1\rangle.\langle
  t^1,t^2\rangle. \cdots .\langle t^{n-1},t^n\rangle$ is a product of
  members of $F$ of the form $s_{\epsilon_1}\circ \cdots \circ
  s_{\epsilon_k}(\langle p_m,q_m\rangle)$ (or its inverse in the
  ``vice versa'' case) and hence $\langle r,s\rangle$, which is an
  arbitrary element of Assoc$(S)$, lies in $H$; i.e.\ Assoc$(S)\leq
  H$.  Since Assoc$(S)$ contains $K$ and (as we have already said) is
  invariant under both shifts, we also have $H\leq \text{Assoc}(S)$.
\end{proof}

In Example~\ref{exam:not normal} we saw a non-normal \Assoc($S$),
expressed as the stabilizer of a subset of $I$. There are of course
many other such non-normal subgroups, for example many expressible as
stabilizers.

\section{Solvable bracket algebras}

We say that the bracket algebra $S$ is \defem{solvable} if there exist
$0\in S$, $n\geq 1$, and $u\in\Tau_n$ such that
$\widehat{u}(S^n)=\{0\}$.  Note that if $v\in\Tau_m$ is an expansion
of $u$ then $\widehat{v}(S^m)=\{0\}$.

\begin{theorem}\label{theorem:solvable} 
  Let $S$ be a solvable bracket algebra.  Any strongly regular law
  $p\ident q$ eventually holds in $S$.
\end{theorem}

\begin{proof}   
  We write $r\geq s$ if the binary tree $r$ is an expansion of the
  binary tree $s$.  This is a partial ordering.  Find expansions
  $\beta_1$, $\beta_2$ such that $\beta_1(p)\geq u$, and
  $\beta_2(\beta_1(q))\geq u$.  Let $\beta=\beta_2\circ\beta_1$.  Then
  $\beta(p)=\beta_2(\beta_1(p)\geq\beta_1(p)\geq u$, and
  $\beta(q)=\beta_2(\beta_1(q))\geq u$.  Since $\beta(p),\beta(q)\geq
  u$ we have $\widehat{\beta(p)}=0=\widehat{\beta(q)}$ and $S\models_e
  p\ident q$.
\end{proof}

This theorem implies that if $S$ is solvable then \Assoc($S$) $=F$.
In the case of groups or Lie algebras our definition of ``solvable" is
easily seen to be equivalent to the usual definitions of ``solvable"
in those contexts.  For solvable Lie algebras and solvable groups the
group of stable associativities is $F$.

The proofs of Propositions~\ref{prop:groups statement}
and~\ref{prop:general stuff} are now complete.

\section{Lie Algebras}

In this section we prove Theorem~\ref{thm:fd lie algebras statement}
and Addendum~\ref{Kac-Moody addendum}.

\begin{lemma}\label{lemma:pre sl2}
  The four element bracket algebra $S=\{0,a,b,c\}$ with operation
  table
  \[
  \begin{array}{c|cccc}
      & 0 & a & b & c \\ \hline
    0 & 0 & 0 & 0 & 0 \\
    a & 0 & 0 & a & b \\
    b & 0 & a & 0 & c \\
    c & 0 & b & c & 0 \\
  \end{array}
  \]
  satisfies no non-trivial strongly regular law, so \Assoc($S$) is trivial.
\end{lemma}

\begin{proof}
  We note two things from this table. First, the function
  $\widehat{[,]}:S\times S\to S$ is surjective; i.e. $S$ is simply
  perfect. Therefore for any $p\in\Tau$ the $n$-ary function
  $\widehat{p}$ is surjective.  Secondly, for any two non-zero $u\neq
  v\in S$, the centralizer of $\{u,v\}$, $C_S(\{u,v\})=\{x\in
  S|[x,u]=[x,v]=0\}$ is trivial, i.e. $\{0\}$.

  Fix $p\in\Tau_n$. Given $u\in S$ we denote by $\widehat{p}|_{x_i=u}$
  the restriction of $\widehat{p}$ to the subset
  $S^{i-1}\times\{u\}\times S^{n-i}\subseteq S^n$; we think of it as
  fixing the value of one of the arguments of $\widehat{p}$.
  Similarly for $\widehat{p} |_{x_i=u,x_j=v}$ where $1\leq i \neq j
  \leq n$ and $u,v\in S$.  We will denote the image of
  $\widehat{p}|_{x_i=u}$ by $\widetilde{p}|_{x_i=u}$. Let $d_i$ denote
  the depth of the $i$-th leave of $p$.

  {\it Claim 1}: If $d_i=1$ then for any $u\in S-\{0\}$ the set
  $\widetilde{p}|_{x_i=u}$ is non-trivial (contains 2 of the values in
  $S-\{0\}$) and has trivial centralizer.  It follows that if $d_i\geq
  1$ then $\widetilde{p}|_{x_i=u}$ is non-trivial and has trivial
  centralizer.

  {\it Claim 2}: If the $i$-th and $(i+1)$-th leaves do not form a
  free caret, i.e. $l_i$ and $l_{i+1}$ differ at more than the last
  bit, then for any $u,v\in S-\{0\}$ the function
  $\widehat{p}|_{x_i=u,x_{i+1}=v}$ is not trivial.  To see this, WLOG
  assume that the $i$-th leaf is in the left subtree $p_0$ of $p$, and
  the $({i+1})$-th leaf is in the right subtree $p_1$ of $p$,
  i.e. that $l_i$ and $l_{i+1}$ have no non-empty common prefix.
  (Otherwise, apply the result to $p_w$, the subtree of $p$ at the
  vertex $w$, where $w$ is the largest common prefix of $l_i$ and
  $l_{i+1}$ i.e. to the smallest subtree of $p$ that contains both the
  $i$-th and the $({i+1})$-th leaves, and use Claim 1.)  Also WLOG
  assume $d_i\geq d_{i+1}$. Since $l_i$ and $l_{i+1}$ differ at more
  than the last bit it cannot be the case that $d_i=d_{i+1}=1$, so
  $d_i \geq 2$.  Note that the depths of leaves in $p_0$ and $p_1$ are
  one less than the corresponding depths in $p$.  From Claim 1 we know
  that $\widetilde{p_0}|_{x_i=u}$ has trivial centralizer, and
  $\widetilde{p_1}|_{x_1=v}$ is non-trivial.  So
  $\widehat{p}|_{x_i=u,x_{i+1}=v}$ is non-trivial.

  Now fix $p,q\in\Tau_n$. We want to show that if $S\models p\ident q$
  then $p=q$.  Since $S$ is simply perfect, by Remark~\ref{rem:simply
    perfect} we may assume that the pair $(p,q)$ is the reduced
  representative of $\langle p,q\rangle$ in $F$.  Suppose $p\neq q$.
  The non-trivial binary tree $q$ has at least one free caret, say
  with the $i-$th and $(i+1)$-th leaves.  Then the $i-$th and
  $(i+1)$-th leaves of $p$ do not form a free caret.  Select $u\in
  S-\{0\}$.  Then $\widehat{q}|_{x_i=u,x_{i+1}=u}=0$ but
  $\widehat{p}|_{x_i=u,x_{i+1}=u}\neq 0$. So, $S\not\models p\ident
  q$.
\end{proof}

\begin{proposition}\label{lemma:sl2}
  If $L$ is the simple Lie algebra $sl_2$ then \Assoc($L$) is trivial.
\end{proposition}

\begin{proof}
  The bracket operation of $sl_2$ is given by
  \[
  \begin{array}{c|ccc}
           & e_{-1} & e_0     & e_1    \\ \hline
    e_{-1} & 0      & -e_{-1} & -2e_0  \\
    e_0    & e_{-1} & 0       & -e_1   \\
    e_1    & 2e_0   & e_1     & 0      \\
  \end{array}
  \]
  It is clear that the algebra $S$ in Lemma~\ref{lemma:pre sl2} is a
  quotient of a bracket subalgebra of $sl_2$.  Since $S$ satisfies no
  strongly regular law, neither does $sl_2$.
\end{proof}

\begin{proposition}\label{prop:fd lie algebras}
  A finite dimensional complex Lie algebra is either solvable or it
  has a subquotient isomorphic to $sl_2$.
\end{proposition}

\begin{proof}
  If $L$ is not solvable, the solvable radical $R(L)$ is a proper
  ideal of $L$ and $L/R(L)$ is a non-trivial semisimple Lie algebra.
  By Serre's theorem on semisimple Lie algebras, $L/R(L)$ contains a
  copy of $sl_2$.
\end{proof}

We can now prove Theorem~\ref{thm:fd lie algebras statement}.

\begin{proof}
  If $L$ is solvable, use Theorem~\ref{theorem:solvable}.  If $L$ is
  not solvable Proposition~\ref{prop:fd lie algebras} implies that $L$
  has a subquotient isomorphic to $sl_2$.  But Lemma~\ref{lemma:sl2}
  says that the bracket algebra $sl_2$ satisfies no strongly regular
  law; so $L$ doesn't either.
\end{proof}

Addendum~\ref{Kac-Moody addendum} follows from the proof since every
infinite dimensional Kac-Moody algebra contains a copy of $sl_2$.

\section{The groups $A_5$ and $S_{\infty}$}

We first prove Proposition~\ref{prop:A5}.  

\begin{proof}
  Note that $[,]:A_5\times A_5\to A_5$ is surjective.  Fix
  $p\in\Tau_n$.  We will use the same notation as in the proof of
  Lemma~\ref{lemma:pre sl2}.

  Claim 1 of that lemma is replaced by: If $d_i=1$ then for any $1\neq
  a\in A_5$, $\widetilde{p}|_{x_i=a}$ is a set with 12, 15, or 30
  elements, depending on whether $v$ is a $5$-cycle, a $3$-cycle, or a
  product of two disjoint $2$-cycles.  It follows that if $d_i\geq 1$
  then $\widetilde{p}|_{x_i=a}$ has trivial centralizer.  This follows
  from direct calculation. The rest of the argument is the same as in
  Lemma~\ref{lemma:pre sl2}.

\end{proof}

Next we prove Proposition~\ref{proposition:Sinfinity}.  Recall that
$S_\infty$ denotes the group of permutations of $\N$ with finite
support.  This group is not finitely generated but it is obviously
elementary amenable, being the union of finite groups.

\begin{proof}
  It is obvious that $S_\infty$ has trivial group of associativities,
  since it contains a copy of $A_5$.  We will show that any $H\leq G$
  of finite index also contains a copy of $A_5$.

  Since there are infinitely many pairwise disjoint 5-cycles in
  $S_\infty$, there are two in the same coset of $H$, so $H$ contains
  a product of two disjoint 5-cycles.  In fact, $H$ contains
  infinitely many such products, and these infinitely many can be
  chosen to be pairwise disjoint.  Now for each of these products of
  two disjoint 5-cycles in $H$ pick two symbols in one of the 5-cycles
  and consider the transposition of those two symbols. There are
  infinitely many of these transpositions, so two of them are in the
  same coset of $H$, and $H$ contains a product of two such
  transpositions.  So we have elements of $H$ of the form
  \[\begin{array}{rcl}
  x&=&(a_1 a_2 a_3 a_4 a_5)(b_1 b_2 b_3 b_4 b_5) \\
  y&=&(c_1 c_2 c_3 c_4 c_5)(d_1 d_2 d_3 d_4 d_5) \\
  z&=&(a_1 a_2)(c_1 c_2)\\
  \end{array}
  \]

  Considering the action of the subgroup $\langle xy,z\rangle$ on the
  set $\{a_1,a_2,a_3,a_4,a_5\}$ we see that $S_5$ is a quotient of
  this subgroup.  By Propositions~\ref{prop:general stuff}
  and~\ref{prop:A5} it follows that $S_{\infty}$ is large.
\end{proof}

\section{The Five Variable Law}

Recall that the Five Variable Law is represented by the following
picture, and in $F$ by the element $c_0 = [x_0, x_1]$ :

\[ \left( \lczero\ , \rczero \right) \]

We are to prove Theorem~\ref{thm:eventual FVL}, namely that
\Assoc($S$) is a non-trivial normal subgroup of $F$ if and only if
this law holds eventually.  (Recall that by Corollary~\ref{cor:FVL on
  the nose} the word ``eventually" is unnecessary in the case of a
simply perfect bracket algebra.)  The ``only if" part is clear since
every non-trivial normal subgroup of $F$ contains the commutator
subgroup.  We prove ``if".

Again we use the PL homeomorphism model of $F$.  It is well known that
the commutator subgroup $F'$ consists of those homeomorphisms in $F$
which agree with the identity map near 0 and near 1.  Let $F_k$ denote
the subgroup of $F$ consisting of homeomorphisms supported on the
closed interval $[\frac {1}{2^k}, 1-\frac {1}{2^k}]$ and let $F_k^+$
denote the subgroup generated by $F_k\cup s_{0}(F_k) \cup s_1(F_k)$.

\begin{lemma}\label{lemma:support} 
  For each $k \geq 2$ we have $F_k^+=F_{k+1}$.
\end{lemma}

Postponing the proof, this lemma implies that the smallest subgroup of
$F$ containing $F_2$ and invariant under the shifts $s_0$ and $s_1$ is
the commutator subgroup $F'$.  Now, inspection of pictures shows that
$F_2$ is a copy of $F$ generated by $c_0$ and $c_1$, the latter being
defined by the following picture (it is just $x_1$ concentrated on the
interval $[\frac {1}{4}, \frac {3}{4}]$):

\[ c_1= \left( \lcone\ , \rcone \right) \]

As self-homeomorphisms of the closed unit interval $I$ we have:

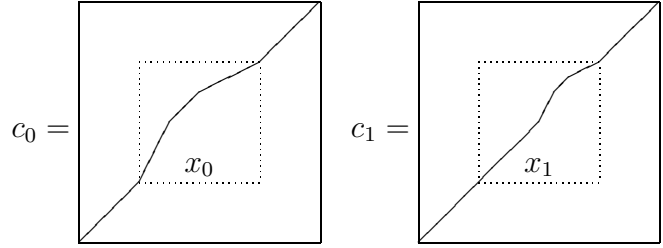
\[ c_0= \raisebox{-40pt}{\plczero} \quad c_1= \raisebox{-40pt}{\plcone} \]

But $c_1$ is just the commutator $[c_0,s_1(c_0)]$, so the smallest
subgroup containing $c_0$ and invariant under both shifts is $F'$.
The ``if" part of Theorem~\ref{thm:eventual FVL} follows.

It remains to prove the Lemma:

\begin{proof}  
  Let $F_{k+1}(\frac{1}{2})$ denote the subgroup of $F_{k+1}$ which
  fixes pointwise the closed interval of length $\frac{1}{2^k}$ whose
  center point is $\frac{1}{2} \in I$.  This is the subgroup generated
  by $s_{0}(F_k) \cup s_1(F_k)$.  Starting with $h \in F_{k+1}$, we
  wish to ``work on" $h$, i.e. to compose $h$ with some members of
  $F_k^+$, until we get a homeomorphism in $F_{k+1}(\frac{1}{2})$;
  that will be enough.  (We sketch the idea since a detailed write-up
  only obscures it; this will prove one of the inclusions and the
  other is obvious.)

  The first step is to follow $h$ by a member of $F_k^+$ to ``move"
  $h(\frac {1}{2})$ to $\frac {1}{2}$. This can be done with a member
  of $F_k$ if $h(\frac{1}{2})$ lies in the open interval $(\frac
  {1}{2^k}, 1-\frac {1}{2^k})$; if that is not so, use an element of
  $s_{0}(F_k) \cup s_1(F_k)$ to make it so.  In summary, we replace
  $h$ by $h_1$ which fixes the point $\frac {1}{2}$.  We then work
  further on $h_1$ without altering it at $\frac {1}{2}$ so as to
  replace it by $h_2$ whose slope is 1 on the closed interval of
  length $\frac {1}{2^{k+1}}$ lying immediately to the left of $\frac
  {1}{2}$; this may involve working on $h_1^{-1}$ instead, but that is
  just as good.  Finally, similar further work gets us to $h_3$ (or
  its inverse) whose slope is also 1 on the closed interval of length
  $\frac {1}{2^{k+1}}$ immediately to the right of $\frac {1}{2}$.
  That is an element of $F_{k+1}(\frac{1}{2})$.
\end{proof}

\section{Bracket algebras having two-sided identity elements}

Here we prove Theorem~\ref{thm:identity element}.  For this a change
of notation is convenient.  Previously, we denoted the depth of the
$i$-th leaf by $d_i$.  In what follows only the depths of the first
(leftmost) and last (rightmost) leaves of a tree $p$ are important;
here we denote those depths by $l(p)$ and $r(p)$.

\begin{proof}
  Let $p,q\in\Tau_n$ be such that $S\models p\ident q$. We want to
  show that $p=q$.  Define $p^l,p^r\in\Tau_{n-1}$ to be such that
  $\widehat{p^l}(x_1,\dots,x_{n-1})=\widehat{p}(1,x_1,\dots,x_{n-1})$
  and
  $\widehat{p^r}(x_1,\dots,x_{n-1})=\widehat{p}(x_1,\dots,x_{n-1},1)$.
  Similarly for $q$.  Then by induction on $n$ we have $p^l=q^l$ and
  $p^r=q^r$. \\

  Case 1: If $l(p)=l(q)=1$ then $p^l=p_1$ and $q^l=q_1$, so we have
  $p_1=q_1$ and $p=q$.  Similarly, if $r(p)=r(q)=1$ then $p=q$.\\

  Case 2: $1<l(p)$ and $1< l(q)$; then $r(p)=r(p^l)=r(q^l)=r(q)$.
  Subcase 2a: if $r(p)=r(q)=1$ then $p=q$ by Case 1.  Subcase 2b: If
  $r(p)=r(q)>1$ and $l(p)=l(q)=1$ then $p=q$ by Case 1.  Subcase 2c:
  If $r(p)=r(q)>1$ and $l(p)=l(q)>1$ then $p_0=(p^r)_0=(q^r)_0=q_0$
  and $p_1=(p^l)_1=(q^l)_1=q_1$, so $p=q$.  Similarly, if $1<r(p)$ and
  $1< r(q)$ then $p=q$.\\

  Case 3: If $l(p)=1< l(q)$ [resp. $r(q)=1<r(p)$] then
  $r(p)=1+r(p_1)=1+r(p^l)=1+r(q^l)=1+r(q)$[resp. $l(q)=1+l(p)$ by
    similar reasoning]. By Case 2, we cannot have $r(p)>1$, so
  $r(p)=1$ and $r(q)=2$[resp. $l(p)=1$ and $l(q)=2$].  Now,
  $p_1=p^l=q^l=[q_{01},\cdot]$, so $p_{10}=q_{01}$, and $p\ident q$ is
  an expansion of the three variable associative law.  But then, since
  $S$ is simply perfect (because it has an identity element), by
  Remark~\ref{rem:simply perfect}, the assertion $S\models p\ident q$
  implies that $S$ is associative, contrary to our hypothesis. So Case
  3 does not arise.\\ 
\end{proof}

A one-sided identity element is not enough in
Theorem~\ref{thm:identity element}.

\begin{example}\label{exam:small not normal}
  Consider the bracket algebra $S(4)=\{1,a,b,c\}$ with binary
  operation $[x,1]=x$, $[x,a]=b$, $[x,b]=c$ and $[x,c]=c$, for any
  $x\in S(4)$. It has a right identity element, $1$. It is not
  associative since $[[x,y],a]=b$, but $[x,[y,a]]=[x,b]=c$.  However,
  it satisfies the law $[x,[[y,z],w]]=[x,[y,[z,w]]]$.  To see this
  consider the cases $w=1$ and $w\neq 1$.  If $w=1$, then both sides
  reduce to $[x,[y,z]]$.  If $w\neq 1$, then both sides equal $c$.
\end{example}

If $(p,q)$ is the usual representative pair of 4-leaf trees for the
standard generator $x_1$ then $\widehat{p}=\widehat{q}$ represents the
law $[x,[[y,z],w]]=[x,[y,[z,w]]]$ which holds in $S(4)$.  Thus $S(4)$
is a quotient of the bracket algebra $S$ in Example~\ref{exam:not
  normal}.  We saw that \Assoc($S$) is a non-normal subgroup of $F$.
We remark that \Assoc($S(4)$) is that same group.  We omit the proof
of this.

Ross Geoghegan, Binghamton University (SUNY), Binghamton  N.Y., 13902,
U.S.A.\\

Fernando Guzm\'an, Binghamton University (SUNY), Binghamton  N.Y., 13902,
U.S.A.\\


\begin{thebibliography}{99}

\bibitem{A} M. Abert, {\em Group laws and free subgroups in topological groups} Bull. London Math. Soc. (to appear).  

\bibitem{BG} K.S.Brown and R.Geoghegan,
{\em An infinite-dimensional torsion free $FP_{\infty}$ group},
Inventiones Math. 77 (1984) 367-381.

\bibitem{BGS}  L. Bartholdi, R. Grigorchuk, and Z.  \v Suni\'k, {\em
Branch groups}.  Handbook of algebra, Vol. 3,  989--1112, North-Holland,
Amsterdam, 2003.


\bibitem{BS} M. Brin and C. C. Squier,  {\em Groups of piecewise linear
homeomorphisms of the real line}, Inventiones. Math. 79 (1985) 485--498.


\bibitem{CFP} J.W. Cannon, W.J. Floyd, and W.R. Parry,{\em Introductory
notes on Richard Thompson's groups}, Enseign. Math. (2) 42 (1996),
215--256.

\bibitem{Gr} R. I. Grigorchuk, {\em On Burnside's problem on periodic
groups} Funktsional. Anal. i Prilozhen.  14 (1980), 53--54.

\bibitem{S} Z. Sunik, private communication.

\end{thebibliography}
\end{document}